\documentclass[12pt]{amsart}
\usepackage{amsfonts}
\usepackage{enumerate}
\usepackage{helvet,courier}
\usepackage{mathptmx,amsmath}
\usepackage{type1cm}
\DeclareMathAlphabet{\mathcal}{OMS}{cmsy}{m}{n}
\DeclareMathAlphabet{\mathbbold}{U}{bbold}{m}{n}  
\usepackage{amsthm}
\usepackage{graphicx}       
\usepackage{multicol}
\usepackage{mathrsfs}
\usepackage{amssymb}
\usepackage{verbatim}
\usepackage{esint}

\usepackage[usenames,dvipsnames]{xcolor}

\usepackage[normalem]{ulem} 

\theoremstyle{plain}
\newtheorem{thm}{Theorem}[section]
\newtheorem{lm}[thm]{Lemma}

\newtheorem{prop}[thm]{Proposition}

\newtheorem{con}{Conjeture}

\theoremstyle{remark}
\newtheorem{rmk}{Remark}

\numberwithin{equation}{section}

\newcommand{\bnu}{\begin{enumerate}}
\newcommand{\enu}{\end{enumerate}}
\newcommand{\bpf}{\begin{proof}}
\newcommand{\epf}{\end{proof}}

\newcommand{\q}{\quad}
\newcommand{\qq}{\qquad}

\newcommand{\sset}{\subset}

\newcommand{\al}{\alpha}
\newcommand{\be}{\beta}
\newcommand{\ga}{\gamma}
\newcommand{\Ga}{\Gamma}

\newcommand{\la}{\lambda}

\newcommand{\ep}{\epsilon}

\newcommand{\vp}{\varphi}
\newcommand{\de}{\delta}

\newcommand{\bbz}{\mathbb{Z}}

\newcommand{\bbr}{\mathbb{R}}

\newcommand{\rn}{{\mathbb{R}^n}}

\newcommand{\bbn}{\mathbb{N}}

\newcommand{\ges}{\gtrsim}
\newcommand{\les}{\lesssim}

\newcommand{\f}{\frac}

\newcommand{\p}{\partial}
\newcommand{\nf}{\infty}

\newcommand{\tf}{\tfrac}
\newcommand{\wh}{\widehat}

\newcommand{\spt}{\text{supp\ }}

\allowdisplaybreaks

\begin{document}
 \author{Peng Chen}
 
  \address{Peng Chen, Department of Mathematics, Sun Yat-sen
 University, Guangzhou, 510275, P.R. China}
 \email{chenpeng3@mail.sysu.edu.cn}
 
\author{Danqing He}
\address{Danqing He, School of Mathematical Sciences,
Fudan University, People's Republic of China}
\email{hedanqing@fudan.edu.cn}

 \author{Xiaochun Li}
 \address{Xiaochun Li, Department of Mathematics, University of Illinois at Urbana-Champaign, Urbana, IL, 61801, USA}
 \email{xcli@math.uiuc.edu}

\author{Lixin Yan}

 \address{Lixin Yan, Department of Mathematics, Sun Yat-sen   University,
 Guangzhou, 510275, P.R. China}
 \email{mcsylx@mail.sysu.edu.cn}

\thanks{
D. He is supported by  National Key R$\&$D Program of China (No. 2021YFA1002500), NNSF of China (No. 12322105), the New Cornerstone Science Foundation, and Natural Science Foundation of Shanghai (No. 22ZR1404900 and 23QA1400300). P. Chen was supported by NNSF of China 12171489 and Guangdong Natural Science Foundation 2022A1515011157. X. Li is partially supported by Simons foundation. L. Yan was  supported by National Key R$\&$D Program of China 2022YFA1005700.
}

\title{On pointwise convergence of cone multipliers}

\date{}

\maketitle

\begin{abstract}
For $p\ge 2$, and $\la>\max\{n|\tf 1p-\tf 12|-\tf12, 0\}$, we prove the pointwise convergence of cone multipliers, i.e. $$ \lim_{t\to\nf}T_t^\la(f)\to f \text{  a.e.},$$ where  $f\in L^p(\rn)$ satisfies $\spt\wh f\sset\{\xi\in\rn:\ 1<|\xi_n|<2\}$. Our main tools are weighted estimates for maximal cone operators, which are consequences of trace inequalities for cones.

\end{abstract}

\section{Introduction}

The Bochner-Riesz operator  is defined by
$$
B^\la(f)(x)=\int_\rn (1-|\xi|^2)^\la_+\wh f(\xi)e^{2\pi ix\cdot\xi} d\xi
$$
for $\la> 0$.
It has been  conjectured that $B^\la$ is bounded on $L^p(\rn)$ when
$\la>\max\{|\tf np-\tf n2|-\tf12,0\}$.
The conjecture was resolved when $n=2$ by  \cite{Carleson1972}; later \cite{Hoermander1973} gave another proof.
 Due to its significant connection with the restriction theorem and other conjectures in harmonic analysis, the Bochner-Riesz conjecture has been extensively studied over the past decades;  see  \cite{Bochner1936}, \cite{Fefferman1970}, \cite{Fefferman1971}, \cite{Fefferman1973a},
 \cite{Bourgain1991a}, \cite{Stein1993}, \cite{Lee2004}, \cite{Guo2021}
and references therein.

A cone multiplier in $\rn$  for $\la> 0$ is defined as
\begin{equation}\label{e041}
\tilde m^\la (\xi)=(1-\tf{|\xi'|^2}{\xi_n^2})^\la_+,
\end{equation}
where $\xi=(\xi',\xi_n)\in\rn$, and the related operator is defined by 
$$
\tilde T^\la(f)(x):=\int_{\rn} \tilde m^\la(\xi)\wh f(\xi) e^{2\pi ix\cdot\xi}d\xi.
$$
It was Stein \cite{Stein1979} who introduced this 
operator.
He conjectured that the range of boundedness 
of $\tilde T^\la$ in $\rn$ coincides with that of $B^\la$ is $\bbr^{n-1}$, as suggested by de Leeuw's theorem  \cite{Leeuw1965}.

On the study of Stein's conjecture in $\bbr^3$, Mockenhaupt
\cite{Mockenhaupt1993a} reduced it to a reverse square function estimate, 
initially obtaining a non-optimal result using a geometric observation on the distribution of plates.
 Bourgain
\cite{Bourgain1995} later improved   Mockenhaupt's result by introducing the so-called bilinear method. This method was  systematically developed to study the restriction theorem and related conjectures;  for further details, see 
\cite{Tao2000}, \cite{Tao2000a}, \cite{Bennett2006}, \cite{Bourgain2011}, \cite{Lee2012b},  \cite{Guth2016}, and \cite{Ou2017}.
The  reverse square function estimate also  implies  the well-known  local smoothing estimates; see \cite{Mockenhaupt1992}. 
To tackle these problems,  Wolff \cite{Wolff2000} introduced decoupling inequalities, which have become a crucial topic in modern harmonic analysis, appearing naturally in various fields from PDEs to number theory; see \cite{Laba2002}, \cite{Garrigos2009}, \cite{Heo2011}, \cite{Bourgain2014a}, \cite{Bourgain2015}, \cite{Du2016}, \cite{Du2018a}, \cite{Guth2018a}, \cite{Guo2019} and references therein. The sharp reverse square function estimate in $\bbr^3$ was recently obtained by Guth, Wang and Zhang 
\cite{Guth2019a}, thereby resolving the cone multiplier conjecture and the local smoothing conjecture in  
 $\bbr^3$, which was generalized to the variable coefficients setting in \cite{Gao2023}.
More results and methods on the boundedness of cone multipliers can be found in \cite{Hong2000},  \cite{Heo2009}, and  \cite{Heo2010}.\\

In this exposition, we shift our focus to the pointwise convergence of cone multipliers, which is connected to maximal estimates of cone multipliers. Pointwise convergence of Fourier series and Fourier integrals is a central issue in harmonic analysis. For classical results, one may refer to  \cite{Zygmund2002}.  According to Stein's maximal principle  \cite{Stein1961a}, 
pointwise convergence problems are consequences of the boundedness problems of corresponding maximal operators. For instance, Carleson  \cite{Carleson1966}  demonstrated the boundedness of the Carleson operator, resolving the long-standing Lusin's conjecture regarding the pointwise convergence of Dirichlet means in $\bbr$; a different proof was given by Fefferman \cite{Fefferman1973}.  
Carbery \cite{Carbery1983} established the $L^4(\bbr^2)$ boundedness of maximal Bochner-Riesz operator. 
For recent progress, interested readers are directed to \cite{Tao2002}, \cite{Lee2004}, and \cite{Li2019}.\\

We formulate our problem as follows. 
Since $\tilde m^\la(t^{-1}\xi)=\tilde m^\la(\xi)$, we dilate $\tilde m^\la$ anisotropically by considering $\tilde m(t^{-1}\xi',\xi_n)$ for $t>0$, which converges to 
$\chi_{\rn}(\xi)$ almost everywhere as $t\to \nf$.  This observation, combined with the Lebesgue dominated convergence theorem,  implies that 
\begin{equation}\label{e0112}
\tilde T_t^\la(f)(x):=\int_{\rn} \tilde m^\la(t^{-1}\xi',\xi_n)\wh f(\xi) e^{2\pi ix\cdot\xi}d\xi
\end{equation}
converges to $f(x)$ for all $x\in\rn$ as $t\to\nf$ when $f$ is a Schwartz function.
It is natural to inquire about the sharp range of  $p$ such that 
\begin{equation}\label{e035}
\lim_{t\to\nf} \tilde T_t^\la(f)(x)=f(x) \qq \text{a.e.} \qq \forall f\in L^p(\rn).
\end{equation}
By Stein's maximal principle, this result can be reduced to the study  of the maximal cone operator
$$
\tilde T_*^\la(f)(x):=\sup_{t>0}|\tilde T_t^\la(f)(x)|.
$$

Since  $\tilde T^\la(f)\le \tilde T_*^\la(f)$, concerning the pointwise convergence and maximal estimates of cone multipliers, 
we pose the following conjecture.

\begin{con}
Let $n\ge 3$, and $p\ge 2$. The  range of $p$ such that \eqref{e035} holds is the same as  the range of $p$ such that  $\|\tilde T^\la\|_{L^p(\rn)\to L^p(\rn)}<\nf$, 
i.e. $\la>\max\{(n-1)|\tf 1p-\tf 12|-\tf12, 0\}$.   
Moreover, $\|\tilde T^\la_*\|_{L^p(\rn)\to L^p(\rn)}<\nf$ when
 $\la>\max\{(n-1)|\tf 1p-\tf 12|-\tf12, 0\}$.

\end{con}

A corresponding conjecture for $p<2$ could also be proposed, but we will not delve into it here as the boundedness of the maximal Bochner-Riesz operator is not clear even in the plane.  For the pointwise convergence of $\tilde T_t^\la$, we derive the following restricted result.

\begin{thm}\label{t013}
Let $n\ge 3$, $p\in[2,\nf)$, and $\la>\max\{|\tf np-\tf n2|-\tf12, 0\}$.
Suppose that $f\in L^p(\rn)$ satisfies that 
\begin{equation}\label{e036}
\spt\wh f\sset\{\xi\in\rn:\ a<|\xi_n|<b\}
\end{equation}
for some $0<a<b<\nf$.
Then 
$\tilde T^\la_t(f)\to f$ a.e. as $t\to\nf$.

\end{thm}

\begin{rmk}
This result is non-optimal compared with the conjectured range  $\la>\max\{(n-1)|\tf 1p-\tf 12|-\tf12, 0\}$.
In particular, when $\la$ is greater than the conjectured critical index $(n-2)/2$, whether the pointwise convergence holds for all $p\in[2,\nf)$ is still open.

\end{rmk}

Theorem~\ref{t013} follows obviously from the case $a=1$ and $b=2$.
Therefore, we introduce the following notation:
$$
m^\la(\xi):=(1-\tf{|\xi'|^2}{\xi_n^2})^\la_+\psi(2^{-1}\xi_n),
$$
where $\psi$ is a nonnegative smooth function supported in $[\tf14,1]$ such that
$$
\sum_{\ga\in\bbz}\psi(2^\ga t)=1
$$ 
for $t> 0$. We can   define $T^\la_t$ and $T^\la_*$ correspondingly.
For $f\in\mathcal S(\rn)$ we define 
$$
 T_t^\la(f)(x):=\int_{\rn}  m^\la(t^{-1}\xi',\xi_n)\wh f(\xi) e^{2\pi ix\cdot\xi}d\xi,
$$
$$
 T_*^\la(f)(x):=\sup_{t>0}| T_t^\la(f)(x)|,
$$
and 
$$
L(f)(x)=\int_{\rn} \psi(2^{-1}\xi_n)\wh f(\xi)e^{2\pi ix\cdot\xi} d\xi.
$$

We can actually prove the following result.

\begin{thm}\label{01276}
Let $n\ge 3$, $p\in[2,\nf)$, and $\la>\max\{|\tf np-\tf n2|-\tf12, 0\}$.
Then, for any $f\in L^p(\rn)$, 
$T^\la_t(f)\to Lf$ a.e. as $t\to\nf$.
\end{thm}

Instead of focusing on the $L^p$ boundedness of maximal Bochner-Riesz operators,
Carbery, Rubio de Francia, and Vega \cite{Carbery1988} demonstrated the pointwise convergence of Bochner-Riesz means in all higher dimensions when $p\ge 2$ by establishing weighted $L^2$  estimates of maximal Bochner-Riesz operators. This approach proves to be effective in various  problems; for instance,  Chen, Duong, He, Lee, and Yan \cite{Chen2020} adopt this strategy to establish the pointwise convergence of Bochner-Riesz means for Hermit operators.  We will employ it  to demonstrate that $T_*^\la$ is bounded on some appropriately weighted $L^2$ spaces. 
Specifically, due to the anisotropic dilation in the definition of $T^\la_*$, we utilize the homogeneous weights given by 
\begin{equation}\label{e040}
w_{\al,\be}(x):=|x'|^{-\al}|x_n|^{-\be}
\end{equation}
where 
$x=(x',x_n)\in\bbr^{n-1}\times \bbr$. 
However, these weights introduce some challenging technical complexities. For example, we need to make a careful decompose wherein an $L^p$ function is dissembled into four pieces in Proposition~\ref{116},  and the associated trace inequalities become more intricate. Nevertheless, we can establish the following result.

\begin{thm}\label{06011}
Let $n\ge 3$, $\al\in[0,n-1)$, $\be\in[0,1)$, and $\la>\max\{\tf {\al+\be-1}2,0\}$. 
Then 
\begin{equation}\label{e033}
\|T^\la_*(f)\|_{L^2(w_{\al,\be})}\le C\|f\|_{L^2(w_{\al,\be})}
\end{equation}
for any Schwartz function $f$.
As a result, for any $f\in L^2(w_{\al,\be})$,
we have $T^\la_t(f)\to Lf$ a.e. as $t\to\nf$.

\end{thm}

This result constitutes the primary technical contribution of this study. Since the trace inequalities represent weighted restriction estimates, which play a crucial role in establishing weighted estimates of Bochner-Riesz means as demonstrated in \cite{Carbery1988}, we aim to transform the inequality \eqref{e033} into analogous trace inequalities of cones. However, this transformation necessitates a nuanced argument that relies on the spatial distribution of $f$, owing to the anisotropic dilation inherent in the definition of maximal cone operators; please refer to Section~\ref{t050} for a detailed explanation. To establish the desired trace inequalities, we make use of geometric properties of cones with the help of an argument inspired by the work in \cite{Carbery1988}.

This paper is organized as follows.
In Section~\ref{t081}, we decompose an arbitrary $L^p$ function into $L^2(w_{\al,\be})$ pieces, where we have to keep the Fourier support properties for technical reasons. Section~\ref{t082} discusses the reduction to square functions. In Section~\ref{040} and Section~\ref{t050}, we establish necessary trace inequalities for cones, which lead to desired weighted inequalities for $T^\la_*$ via a delicate argument in term of the physical position of $f$. The proof of main theorems are given in Section~\ref{t083}. 

\medskip

\noindent
{\bf Notations}

By $A\les B$ we mean that there exists an absolute constant $C$ such that $A\le CB$. $A\sim B$ if $A\les B$ and $B\les A$.

$E+O(\de):=\{x\in \rn: dist(x,E)<C\de\}\}$ is the $C\de$-neighborhood of $E$ in $\rn$, where $C$ is an unimportant absolute constant.

$Q(x,\ell)$ is the cube centered at $x\in \rn$ with length $\ell$. We denote  by $NQ(x,\ell)$ the cube $Q(x, N\ell)$.

Given a multiplier $m\in L^\nf(\rn)$, the associated linear operator $T_m$ is defined as
$$
T_m(f)(x)=\int_{\rn} m(\xi)\wh f(\xi)e^{2\pi ix\cdot\xi}d\xi.
$$

$\mathcal S(\bbr^n)$ is the class of Schwartz functions.

\section{Decomposing  functions}\label{t081}

We decompose an arbitrary function $f\in L^p(\rn)$ in terms of the weights $w_{\al,\be}$ defined in \eqref{e040} with the Fourier support remained, therefore, to obtain the pointwise convergence of cone multiplies, it is natural to consider weighted $L^2$ estimates of $T^\la_*$. The proof is standard, but we include it for the sake of completeness.

\begin{prop}\label{116}
Let $p\ge 2$, $\al_1=\be_1=0$, $0\le \al_3<(n-1)(1-\tf2p)<\al_2,\al_4<1+(n-1)(1-\tf2p)$, and $0\le \be_2<1-\tf2p<\be_3,\be_4<2-\tf2p$.
For any $f\in L^p(\rn)$, we can decompose it into
$$
f=f_1+f_2+f_3+f_4
$$ 
such that
$f_i\in L^2(w_i)$,
where $w_i(x)=|x'|^{-\al_i}|x_n|^{-\be_i}$ for $i=1,2,3,4$.
Moreover, if $\spt \wh f\sset S=\{\xi\in\rn:\ \xi_n\in(\tf1{2},2)\}$, then
\begin{equation}\label{e015}
\spt \wh f_i\sset S'=\{\xi\in\rn:\ \xi_n\in(\tf1{10},10)\}
\end{equation} 
for $i=1,2,3,4.$

\end{prop}

\bpf

We take a function $\vp_k\in\mathcal S(\bbr^k)$ such that $\spt \wh\vp_k\sset B(0,10^{-10})$ and $\vp_k(0)=1$.

We define 
$$
f_1(x)=f(x)\vp_{n-1}(x')\vp_1(x_n)
$$
$$
f_2(x)=f(x)(1-\vp_{n-1})(x')\vp_1(x_n)
$$
$$
f_3(x)=f(x)\vp_{n-1}(x')(1-\vp_1)(x_n)
$$
$$
f_4(x)=f(x)(1-\vp_{n-1})(x')(1-\vp_1)(x_n).
$$
Obviously $f=f_1+f_2+f_3+f_4$, and,  when $\spt \wh f\sset S$,
\eqref{e015} can be verified easily using that $\spt \wh\vp_k\sset B(0,10^{-10})$.

Since $p\ge 2$, by H\"older's inequality, we obtain
\begin{align*}
\int_{\rn}|f_1(x)|^2w_1(x)dx= &\int_{\rn}|f(x)\vp_{n-1}(x')\vp_1(x_n)|^2dx\\
\le &C\|f\|_{L^p(\rn)}^2.
\end{align*}

Let us analyze $f_4$ below.
By H\"older's inequality, we control $\|f_4\|_{L^2(w_4)}^2$ by 
\begin{equation}\label{e016}
\|f\|_{L^p}^2\Big(\int_{\rn}((1-\vp_{n-1})(x')|x'|^{-\al_4}(1-\vp_1)(x_n)|x_n|^{-\be_4})^{\tf p{p-2}}dx\Big)^{\tf{p-2}p}.
\end{equation}
As  $\vp_{n-1}(0)=1$, we have 
$$
|1-\vp_{n-1}(x')|\le C|x'|.
$$
This implies
\begin{align*}
&\int_{\bbr^{n-1}}((1-\vp_{n-1})(x')|x'|^{-\al_4})^{p/(p-2)}dx'\\
\le &\int_{|x'|\le 1}|x'|^{(1-\al_4)\tf p{p-2}}dx'+\int_{|x'|\ge 1}|x'|^{-\al_4\tf{p}{p-2}}dx'\\
\les &\int_0^1 r^{(1-\al_4)\tf p{p-2}+n-2}dr+\int_1^\nf r^{n-2-\al_4\tf{p}{p-2}}dr,
\end{align*}
which is finite since $(n-1)(1-\tf2p)<\al_4<1+(n-1)(1-\tf2p)$.
Similarly  
$$\int_{\bbr}((1-\vp_1)(x_n)|x_n|^{-\be_4})^{p/(p-2)}dx_n<\nf
$$
when $\be_4 \in (1-\tf2p,2-\tf2p)$.
These estimates and \eqref{e016} show that
$$
\|f_4\|_{L^2(w_4)}\le C\|f\|_{L^p}.
$$

One can also obtain similarly that
\begin{align*}
\|f_2\|_{L^2(w_2)}^2
=&\int_{\rn}|f_2(x)|^2|x'|^{-\al_2}|x_n|^{-\be_2}dx\\
\le&C\|f\|_{L^p}^2\Big(\int_{\rn}((1-\vp_{n-1})(x')|x'|^{-\al_2}\vp_1(x_n)|x_n|^{-\be_2})^{\tf p {p-2}}dx\Big)^{\tf{p-2}p}\\
<&C\|f\|_{L^p}^2
\end{align*}
when $(n-1)(1-\tf2p)<\al_2<1+(n-1)(1-\tf2p)$ and $\be_2\tf p{p-2}<1$.

The proof for $f_3\in L^2(w_3)$ is similar. 

\epf

\section{Reduction to square functions}\label{t082}

We recall that  $\psi$ is  a smooth bump  supported in $[\tf14,1]$ such that
$\sum_{\ga\in\bbz}\psi(2^\ga t)=1$ for $t> 0$,  
$$
m^\la(\xi)=(1-\tf{|\xi'|^2}{\xi_n^2})^\la_+\psi(2^{-1}\xi_n),
$$
$$
T_t^\la(f)(x)=\int_{\rn} m^\la(t^{-1}\xi',\xi_n)\wh f(\xi) e^{2\pi ix\cdot\xi}d\xi,
$$
and 
$$
T^\la_*(f)(x)=\sup_{t>0}|T^\la_t(f)(x)|.
$$
For  given $\al,\be\ge 0$ and $w_{\al,\be}(x)=|x'|^{-\al}|x_n|^{-\be}$, we  will 
explore the range of $\lambda$ such that 
$T^\lambda_*$ is bounded on $L^2(w_{\al,\be})$.

We define, for $\ga\ge 1$,
$$
m_\ga(\xi)=\psi(2^{\ga}(1-\tf{|\xi'|^2}{\xi_n^2}))m^\la(\xi)2^{\ga\la},
$$
and 
$$
m_0(\xi)=m^\la(\xi)(1-\sum_{\ga\ge 1}\psi(2^\ga(1-\tf{|\xi'|^2}{\xi_n^2}))).$$
We remark that $m_0$ is 
supported in 
$$
\{\xi\in\rn:\ |\xi'|/\xi_n\in[0,1],\ \xi_n\in[\tf12,2]\},
$$
and $m_\ga$ is    supported in 
$$\{\xi\in\rn:\ |\xi'|/\xi_n\in[1-2^{-\ga},1-2^{-\ga-1}],\ \xi_n\in[\tf12,2]\}
$$ 
for $\ga\ge 1$.
Obviously
\begin{equation}\label{e300}
T^\la_*(f)(x)\le\sum_{\ga\ge0}2^{-\ga\la}M_\ga(f)(x),
\end{equation}
where 
$$
M_\ga(f)(x)=\sup_{t>0}\Big|\int_{\rn}m_\ga(t^{-1}\xi',\xi_n)\wh f(\xi)e^{2\pi ix\cdot\xi}d\xi\Big|.
$$

Concerning the weighted boundedness of $M_\ga$, we can easily prove a non-optimal result.

\begin{prop}\label{111}
For $\ga\ge0$, $\al\in(-(n-1),(n-1))$, and $\be\in (-1,1)$,
\begin{equation}\label{e043}
\|M_\ga(f)\|_{L^2(w_{\al,\be})}\le C_\ga\|f\|_{L^2(w_{\al,\be})}.
\end{equation}
\end{prop}

To prove this proposition, we recall first some definitions and results.

We denote by $Q_k$ a cube in $\bbr^k$, and $M_k$ the standard Hardy-Littlewood maximal function in $\bbr^k$.
A strong maximal function (on $\bbr^{n-1}\times \bbr$):
$$
M_S(f)(x)=\sup_{R\ni x}\f1{|R|}\int_{R}|f(y)|dy,
$$
where the supremum runs over all $R$ of the form $Q_{n-1}\times Q_1$. It is well known that $M_S$ is bounded on $L^p(\bbr^n)$ for $1<p<\nf$
since
\begin{equation}\label{e042}
M_S(f)(x',x_n)\le (M_1\circ M_{k-1})(f)(x',x_n).
\end{equation}

A weight $w$ belongs to the (product) $A_p(\bbr^{n-1}\times \bbr)$  if there exists a finite constant $C_w$ such that for all $R=Q_{n-1}\times Q_1$,
\begin{equation}\label{e047}
\Big(\tf1{|R|}\int_Rw(x',x_n)dx'dx_n\Big)\Big(\tf1{|R|} \int_Rw(x',x_n)^{-1/(p-1)}dx'dx_n\Big)^{p-1}\le C_{w}.
\end{equation}
One can easily verify that $w_{\al,\be}$ is an $A_2(\bbr^{n-1}\times \bbr)$ weight if and only if $\al\in(-(n-1),n-1)$ and $\be\in(-1,1)$. 
Since $M_k$ is bounded on $L^2(w)$ when $w$ is an $A_2(\bbr^k)$ weight,
it follows from \eqref{e042} that 
\begin{equation}\label{e044}
\|M_S(f)\|_{L^2(w_{\al,\be})}\le C\|f\|_{L^2(w_{\al,\be})}
\end{equation}
when $\al\in(-(n-1),n-1)$ and $\be\in(-1,1)$.

\bpf[Proof of Propositio~\ref{111}]
As $m_\ga(\xi)$ is a compactly supported smooth function, it is routine to verify  that 
$$
M_\ga(f)(x)\le C_\ga\ M_{S}(f)(x),
$$
which implies further \eqref{e043} by \eqref{e044}.

\epf

The boundedness of $T^\la_*$ on $L^2(w_{\al,\be})$ is now reduced to estimating the norm  $\|M_\ga\|_{L^2(w_{\al,\be})\to L^2(w_{\al,\be})}$ as $\ga\to\nf$, since we have proved that  $M_\ga$ is bounded on $L^2(w_{\al,\be})$ for all $\ga\ge 0$.

We will focus on the case  $\ga\ge 1$ from now on.

By the Sobolev embedding, we obtain further that 
$$
M_\ga(f)(x)^2\le 2^\ga \mathcal G_\ga(f)(x)\tilde {\mathcal G}_\ga(f)(x),
$$
where 
$$
{\mathcal G}_\ga(f)(x):=\Big(\int_0^\nf\Big|\int_{\rn}m_\ga(t^{-1}\xi',\xi_n)\wh f(\xi)e^{2\pi ix\cdot\xi}d\xi\Big|^2\f{dt}t\Big)^{1/2},
$$
and 
$$
\tilde {\mathcal G}_\ga(f)(x):=\Big(\int_0^\nf\Big|\int_{\rn}\tilde m_\ga(t^{-1}\xi',\xi_n)\wh f(\xi)e^{2\pi ix\cdot\xi}d\xi\Big|^2\f{dt}t\Big)^{1/2}.
$$ 
Here
$\tilde m_\ga(\xi',\xi_n)=(2^{-\ga}\xi')\cdot(\nabla_{n-1}m_\ga)(\xi',\xi_n)$, and
$\nabla_{n-1} m_\ga(\cdot,\xi_n)$ is the gradient of $m_\ga(\cdot,\xi_n)$ as a function defined on $\bbr^{n-1}$.
In summary we obtain 
\begin{equation}\label{e01272'}
T^\la_*(f)(x)\le M_0(f)(x)+\sum_{\ga\ge1}2^{-\ga\la}2^{\ga/2}[{\mathcal G}_\ga(f)(x)\tilde {\mathcal G}_\ga(f)(x)]^{1/2}.
\end{equation}

For $0<\de\ll1$, we fix a smooth function
$\mu_\de(\xi')$ such that 
\begin{equation}\label{e008}
\spt\mu_\de\sset \{\xi'\in\bbr^{n-1}:\ |\xi'|\in[1-\de,1] \}
\end{equation}
and
\begin{equation}\label{e046}
|\partial^{l}\mu_\de(\xi')|\le C_l\de^{-|l|}
\end{equation}
for any $l=(l_1,\dots,l_{n-1})\in\bbn^{n-1}$.
Associated with $\mu_\de$, we define the multiplier
$$
m_\de(\xi',\xi_n)=\mu_\de(\xi'/\xi_n)\psi(2^{-1}\xi_n),
$$
and 
$$
{\mathbf G}_\de(f)(x):=\left(\int_0^\nf\Big|\int_{\rn}m_\de(t^{-1}\xi',\xi_n)\wh f(\xi)e^{2\pi ix\cdot\xi}d\xi\Big|^2\f{dt}t\right)^{1/2}.
$$
As the multipliers $ m_\ga $ and $\tilde m_\ga$ are multipliers of this form with $\de=2^{-\ga}$,
the problem is now reduced to studying the boundedness of  
${\mathbf G}_\de$.

\hfill

Let us examine a simple case first.

\begin{prop}\label{117}
Let $\la>0$.
Then 
\begin{equation}\label{e045}
\|{\mathbf G}_\de(f)\|_{L^2(\rn)}\le \de^{1/2}\|f\|_{L^2(\rn)}
\end{equation}
and 
$T^\la_*$ is bounded on $L^2(\rn).$
\end{prop}

\bpf

By Plancherel's identity, $\|{\mathbf G}_\de(f)\|_{L^2(\rn)}^2$  is equal to
$$
\int_\rn\int_0^\nf |m_\de(t^{-1}\xi',\xi_n)|^2 \tf{dt}t |\wh f(\xi)|^2d\xi.
$$
For a fixed $\xi$, the $t$ such that $m_\de(t^{-1}\xi',\xi_n)\neq 0$ is contained in $\tf{|\xi'|}{\xi_n}+O(\tf{|\xi'|}{\xi_n}\de)$, hence the inner integral is bounded by $C\de$, which implies \eqref{e045}.

Recalling \eqref{e01272'} and Proposition~\ref{111}, we obtain from \eqref{e045} that 
$$\|T^\la_*(f)\|_{L^2(\rn)}\le\sum_{\ga\ge 0}2^{-\ga \la}2^{\ga/2}2^{-\ga /2}\|f\|_{L^2(\rn)}\le C\|f\|_{L^2(\bbr^n)}$$
when $\la>0$.
\epf

We  study the boundedness of ${\mathbf G}_\de$ on $L^2(w_{\al,\be})$ below.
Let us fix $\de\ll1$, and decompose ${\mathbf G}_\de$ further. 
Let $S_t(f)$ be the linear operator with multiplier $m_\de(t^{-1}\xi',\xi_n)$, then we can write
$\big[{\mathbf G}_\de(f)(x)\big]^2$ as $\sum_{k\in\bbz}|G_k(f)(x)|^2$ with
$$
|G_k(f)(x)|^2:=\int_{2^k}^{2^{k+1}}|S_t(f)(x)|^2\f{dt}t.
$$

The study of  ${\mathbf G}_\de$ can be reduced to the boundedness of $G_k$, as we will see in the next result.

\begin{prop}\label{101}
Let $\al\in(-(n-1),n-1)$ and $\be\in(-1,1)$. 
If 
\begin{equation}\label{e048}
\|G_0(f)\|_{L^2(w_{\al,\be})}\le A\|f\|_{L^2(w_{\al,\be})}
\end{equation}
for any $f\in L^2(w_{\al,\be})$, then 
$$\|{\mathbf G}_\de(f)\|_{L^2(w_{\al,\be})}\les A\|f\|_{L^2(w_{\al,\be})}.$$

\end{prop}

To prove this proposition, we establish first the following lemma.

\begin{lm}\label{102}
Fix $\al,\be\in\bbr$.
If 
\begin{equation}\label{e202}
\|G_0(f)\|_{L^2(w_{\al,\be})}\le A\|f\|_{L^2(w_{\al,\be})}
\end{equation}
for any $f\in L^2(w_{\al,\be})$, then for any $k\in\bbz$ it holds
$$\|G_k(f)\|_{L^2(w_{\al,\be})}\le A\|f\|_{L^2(w_{\al,\be})}.$$

\end{lm}

\bpf
We need to estimate
\begin{align*}
\int_\rn|G_k(f)(x)|^2w_{\al,\be}(x)dx=&
\int_\rn\int_{2^k}^{2^{k+1}}\Big|S_t(f)(x)\Big|^2\f{dt}t
|x'|^{-\al}|x_n|^{-\be}dx.
\end{align*}
By 
change of variables,
the inner integral is equal to
\begin{align*}
&\int_1^2\Big|\big[m_\de(2^{-k}t^{-1}\xi',\xi_n)\wh f(\xi)\big]^\vee(x)\Big|^2\f{dt}t\\
=&
\int_1^2\Big|2^{k(n-1)}\big[m_\de(t^{-1}\xi',\xi_n)\wh f(2^k\xi',\xi_n)\big]^\vee(2^kx', x_n)\Big|^2\f{dt}t.
\end{align*}
Let $\wh{f_k}(\xi)=\wh f(2^k\xi',\xi_n)$. It follows from \eqref{e202} that
\begin{align*}
&\int_\rn|G_k(f)(x)|^2w_{\al,\be}(x)dx\\
=&
2^{k(n-1)} 2^{k\al}
\int_\rn\int_1^2\Big| S_t(f_k)(x)\Big|^2\f{dt}t
|x'|^{-\al}|x_n|^{-\be}dx\\
\le &A^22^{k(n-1)} 2^{k\al}
\int_\rn |f_k(x)|^2|x'|^{-\al}|x_n|^{-\be}dx\\
=&A^22^{k(n-1)} 2^{k\al}
\int_\rn |2^{-k(n-1)}f(2^{-k}x',x_n)|^2|x'|^{-\al}|x_n|^{-\be}dx\\
\le& A^2\|f\|_{L^2(w_{\al,\be})}^2.
\end{align*}
\epf

To finish the proof of Proposition~\ref{101}, we need the following  lemma to paste $G_k$.
Let $\wh\zeta$ be a smooth bump supported in $\{\xi'\in\bbr^{n-1}:\ |\xi'|\sim 1\}$, and 
$\wh\psi$ a smooth bump supported in $\{\xi_n\in\bbr:\ \xi_n\sim 1\}$. We define $\wh{f_{k,L}}(\xi):=\wh f(\xi)\zeta(2^{-(k+L)}\xi')\wh\psi(2^{-L}\xi_n)$. We will explore the orthogonality of $\{f_{k,L}\}$ below.

\begin{lm}\label{103}
Let $\al\in(-(n-1),n-1)$ and $\be\in(-1,1)$. Then
$$
\sum_{k\in\bbz}\sum_{L\in\bbz}\int_\rn|f_{k,L}(x)|^2w_{\al,\be}(x)dx\les \int_\rn|f(x)|^2w_{\al,\be}(x)dx.
$$

\end{lm}

\bpf
Let $\wh {F_{k,L}}(\xi)=\wh f(\xi)\wh\zeta(2^{-k}\xi')\wh\psi(2^{-L}\xi_n)$, 
which is different from $\wh {f_{k,L}}(\xi)$ in that they have different locations in $\xi'$.
By a simple change of variables on $(k,L)$, we have
\begin{align*}
&\sum_{k\in\bbz}\sum_{L\in\bbz}\int_\rn |f_{k,L}|^2w_{\al,\be}(x)dx\\
=&\sum_{k\in\bbz}\sum_{L\in\bbz}\int_\rn |F_{k,L}|^2w_{\al,\be}(x)dx\\
\sim&\int_\rn\int_0^1\int_0^1|\sum_k\sum_Lr_k(s_1)r_L(s_2)F_{k,L}(x)|^2ds_1ds_2w_{\al,\be}(x)dx,
\end{align*}
where $r_k(s_1)$ and $r_L(s_2)$ are Rademacher functions, and we use Khintchine's inequality of several variables (see \cite[Appendix C.5]{Grafakos2014b}).

Let us write $\sum_k\sum_Lr_k(s_1)r_L(s_2)F_{k,L}(x)=K_{s_1,s_2}*f$ with 
$$
\wh K_{s_1,s_2}(\xi)=\big(\sum_kr_k(s_1)\wh\zeta(2^{-k}\xi')\big)\big(\sum_Lr_L(s_2)\wh\psi(2^{-L}\xi_n)\big)=:\wh K_{s_1}(\xi')\wh K_{s_2}(\xi_n).
$$
This is a bi-parameter singular integral with associated constants independent of $s_1$ and $s_2$.
Applying \cite[Theorem]{Fefferman1988} we obtain
$$
\int_\rn|K_{s_1,s_2}*f(x)|^2w_{\al,\be}(x)dx
\les \int_\rn |f(x)|^2w_{\al,\be}(x)dx
$$
uniformly in $s_1,s_2\in[0,1]$ as $w_{\al,\be}$ is a product $A_2(\bbr^{n-1}\times \bbr)$ weight in the sense of \eqref{e047}. 
This concludes the proof by integrating over $s_1$ and $s_2$.
\epf

\bpf[Proof of Proposition~\ref{101}]
We write $\|{\mathbf G}_\de(f)\|_{L^2(w_{\al,\be})}^2$ as
$$
\sum_k\int_{2^k}^{2^{k+1}}\int_\rn\Big|[m_\de(t^{-1}\xi',\xi_n)\wh f(\xi)]^\vee(x)\Big|^2w_{\al,\be}(x)dx\f{dt}t,
$$
which, by 
the support of $m_\de$,  is bounded by 
\begin{align*}
&\sum_k\int_{2^k}^{2^{k+1}}\int_\rn\Big|[m_\de(t^{-1}\xi',\xi_n)\wh f_{k,0}(\xi)]^\vee(x)\Big|^2 w_{\al,\be}(x)dx\f{dt}t\\
=&\sum_k\int_\rn\int_{2^k}^{2^{k+1}}\Big|[m_\de(t^{-1}\xi',\xi_n)\wh f_{k,0}(\xi)]^\vee(x)\Big|^2\f{dt}t w_{\al,\be}(x)dx\\
=&\sum_k\int_\rn|G_k(f_{k,0})(x)|^2w_{\al,\be}(x)dx\\
\les& A^2 \sum_k\int_\rn|f_{k,0}(x)|^2w_{\al,\be}(x)dx\\
\les& A^2\int_\rn|f(x)|^2w_{\al,\be}(x)dx.
\end{align*}
In the last two steps we use Lemma~\ref{102} and Lemma~\ref{103} 
\epf

\section{Some trace lemmas }\label{040}

We prove in this section some trace lemmas, whose dual forms, weighted restriction-type estimates, are good substitutes of classical restriction estimates in our setting.

Let $\al\ge0$, $\be\ge0$, and 
$$
\Ga_\de:=\{x\in\rn:\ \text{dist }(x,\Ga)<\de\},
$$
where $\Ga=\{\xi\in\rn:\ |\xi'|=\xi_n\in(1,2)\}$.
We are interested in determining $C_\de(\al,\be)$ in  the following inequality.
\begin{equation}\label{e005}
\int_{\Ga_\de}|\wh f(\xi)|^2\le
C_{\al,\be} C_\de(\al,\be)\int_\rn|f(x)|^2w^{-1}_{\al,\be}(x)dx,
\end{equation}
where $w^{-1}_{\al,\be}=|x'|^\al|x_n|^\be$.

Our main result in this section is as follows.

\begin{prop}\label{109'}
Let $\al\in[0,n-1)$ and $\be\in[0,1)$. Then \eqref{e005} holds for 
$$
C_\de(\al, \be)\sim \left\{\begin{array}{ll} 
\de^{\al+\be},&\q \al+\be<1\\
\de\log\tf1\de,&\q \al+\be=1\\
\de,& \q\al+\be>1.
\end{array}\right.$$
\end{prop}

The proof of Proposition~\ref{109'} is contained in 
the following results, as the case $\al=\be=0$ is trivial.

\begin{lm}\label{110}
If $\be=0$, then \eqref{e005} holds for 
$$
C_\de(\al, 0)\sim \left\{\begin{array}{ll} 
\de^\al,&\q 0<\al<1\\
\de\log\tf1\de,&\q \al=1\\
\de,& \q1<\al<n-1.
\end{array}\right.$$
\end{lm}

This result is a direct consequence of the following classical trace inequalities for spheres.

\begin{lm}[{\cite[Lemma 3]{Carbery1988}}]\label{105}

Let $\al\in(0,n-1)$. It holds
\begin{equation}\label{e011}
\int_{\mathbb S^{n-2}+O(\de)}|\wh g(\xi)|^2d(\xi)\les C_\de(\al,0)\int_{\bbr^{n-1}}|g(x)|^2|x|^\al dx.
\end{equation}

\end{lm}

\bpf[Proof of Lemma~\ref{110}]

Let 
$$
f_{\xi_n}(x')=\int_{\bbr^{n-1}}\wh f(\xi',\xi_n)e^{2\pi ix'\cdot\xi'}d\xi'=\int_\bbr f(x',x_n)e^{-2\pi ix_n\xi_n}dx_n.
$$
Applying Lemma~\ref{105} we obtain
\begin{align*}
\int_{\Ga_\de}|\wh f(\xi)|^2d\xi
\le&\int_1^2\int_{\mathbb S^{n-2}_{\xi_n}+O(\de)}|\wh f(\xi',\xi_n)|^2d\xi'd\xi_n\\
\les &C_\de(\al,0)\int_1^2\int_{\bbr^{n-1}}|f_{\xi_n}(x')|^2|x'|^\al
dx' d\xi_n,
\end{align*}
where $\mathbb S^{n-2}_{\xi_n}$ is the sphere centered at the origin of radius $\xi_n$.
By Plancherel's identity, this is bounded by 
\begin{align*}
&C_\de(\al,0)\int_{\bbr^{n-1}}\int_\bbr\Big|\int_\bbr f(x',x_n)e^{-2\pi ix_n\xi_n}dx_n\Big|^2d\xi_n|x'|^\al dx'\\
=&C_\de(\al,0)\int_\rn|f(x)|^2|x'|^\al dx.
\end{align*}

\epf

\begin{lm}\label{106}
If $\al=0$ and $\be\in(0,1)$, then \eqref{e005} holds for
$$
C_\de(\al, \be)\sim \de^\be.
$$
\end{lm}

\bpf

We want to prove
\begin{equation}
	\int_{\Ga_\de}|\wh f(\xi)|^2d\xi \le C\de^\be \int_\rn|f(x)|^2|x_n|^{\be}dx,
\end{equation}
which by duality  is equivalent to
\begin{equation}\label{e101}
\int_\rn |\wh g(\xi)|^2\f{d\xi}{|\xi_n|^\be}\les \de^\be \int_{\Ga_\de}|g(x)|^2dx
\end{equation}
for any $g\in L^2(\Ga_\de)$.

We first claim that
\begin{equation}\label{e102}
\int_\bbr |\wh F(\xi_n)|^2\f{d\xi_n}{|\xi_n|^\be}\les \de^\be \int_{I_\de}|F(x_n)|^2dx_n
\end{equation}
for any $F\in L^2(I_\de)$, where $I_\de=[a-\de,a+\de]$ is an interval of length $2\de$.

Take $F(x_n)=g(x',x_n)$, which is supported in $(|x'|-\de,|x'|+\de)$ as $x\in\Ga_\de$, then
$$
\int_{\bbr^{n-1}}\int_\bbr|\wh F(\xi_n)|^2\f{d\xi_n}{|\xi_n|^\be}dx'\les\de^\be\int_\rn|g(x)|^2dx
$$
by \eqref{e102}.
The left hand side, by Plancherel's identity, is equal to 
$$
\int_{\rn} |\wh g(\xi)|^2\f{d\xi}{|\xi_n|^\be},
$$
which yields \eqref{e101}.

It remains to verify \eqref{e102}. Actually $\int_\bbr |\wh F(\xi_n)|^2\f{d\xi_n}{|\xi_n|^\be}$ equals
$$
C_\be\int_\bbr\int_\bbr F(x_n)\overline F(y_n)\f1{|x_n-y_n|^{1-\be}}dx_ndy_n
$$
as $\be\in(0,1)$.
It follows from this expression and the Schur test that 
\eqref{e102} is a consequence of
$$
\sup_{x_n\in[a-\de,a+\de]}\int_{a-\de}^{a+\de}\f1{|x_n-y_n|^{1-\be}}dy_n\les\de^\be,
$$
which is straightforward to check using $\be>0$.
\epf

\begin{lm}\label{109}
If $\al\in(0,n-1)$ and $\be\in(0,1)$, then \eqref{e005} holds for 
$$
C_\de(\al, \be)\sim \left\{\begin{array}{ll} 
\de^{\al+\be},&\q \al+\be<1\\
\de\log\tf1\de,&\q \al+\be=1\\
\de,& \q\al+\be>1.
\end{array}\right.$$
\end{lm}

\bpf
By duality \eqref{e005} is equivalent to 
\begin{equation}\label{e103}
\int_\rn |\wh g(\xi)|^2|\xi'|^{-\al}|\xi_n|^{-\be}d\xi\les C_\de(\al,\be) \int_{\Ga_\de}|g(x)|^2dx.
\end{equation}

By an argument similar to the proof of Lemma~\ref{106}, it suffices to show  that 
\begin{equation}
\sup_{x\in\Ga_\de}\int_{\Ga_\de}|x'-y'|^{\al-(n-1)}|x_n-y_n|^{\be-1}dy\les C_\de(\al,\be).
\end{equation}
The left hand side is
\begin{equation}\label{e104}
\sup_{x\in\Ga_\de}\int_{\Ga_\de-x}|z'|^{\al-(n-1)}|z_n|^{\be-1}dz.
\end{equation}
By rotation in the $x'$-plane, we may assume that $x'=(|x'|,0,\dots,0)\in \bbr^{n-1}$ with $|x'|\in[1,2]$.

To analyze the integral in \eqref{e104}, we define
$$S_0=\{z\in\Ga_\de-x:\ |z'|\le\de\},$$ 
$$S_\ell=\{z\in\Ga_\de-x:\ \ell\de\le|z'|\le (\ell+1)\de\}\q \text{for } 1\le\ell\le \de^{-1}/1000\},$$ 
and 
$$S_\nf=\{z\in\Ga_\de-x:\ |z'|\ge 1/1000\}.$$
Obviously $\Ga_\de-x=\big[\cup_{0\le\ell\le\de^{-1}/1000}S_\ell]\cup S_\nf$. $S_0$ is essentially a ball of radius $\sim\de$, and $S_\ell$ is contained in an annular cylinder for $\ell\ge 1$. 

Let us calculate the integrals over $S_0$ and $S_\nf$ first.

We observe that $0\in \Ga_\de-x$, and, for $z\in \Ga_\de-x$, we have $|z_n|\le 2\de$ when $|z'|\le\de$, therefore
\begin{equation}\label{e108}
\int_{S_0}|z'|^{\al-(n-1)}|z_n|^{\be-1}dz
\les\int_{|z'|\les\de}|z'|^{\al-(n-1)}dz'\int_{|z_n|\les\de}|z_n|^{\be-1}dz_n
\les\de^{\al+\be}
\end{equation}
as $\al,\be>0$.

To estimate the integral over $S_\nf$, we denote the intersection of $S_\nf$ and the hyperplane 
$\{y\in\rn:\ y_n=z_n\}$ by  
$$
E_{z_n}:=\{z'\in\bbr^{n-1}:\ \big||z'+x'|-(z_n+x_n)\big|\les\de,\ |z'|\ge 10^{-3}\},
$$
whose $(n-1)$-dimensional volume is bounded by $\de$ as $E_{z_n}$ is contained in an annulus of radius $\sim1$ and width $\sim\de$.  This yields
\begin{equation}\label{e107}
\int_{S_\nf}|z'|^{\al-(n-1)}|z_n|^{\be-1}dz\les \int_{|z_n|\les 1}|z_n|^{\be-1}\int_{E_{z_n}}|z'|^{\al-(n-1)}dz'dz_n\les\de
\end{equation}
since $|z'|\sim 1$ in $E_{z_n}$ and $\be>0$.

We estimate the integral over $S_\ell$ with $1\le \ell\le \de^{-1}/1000$ below, where we have $\ell\de\le|z'|\le (\ell+1)\de$.
It is easy to see that $|z_n|\le (\ell+5)\de$ in $S_\ell$. Fix $|z_n|\le (\ell+5)\de$, then the set 
$\{z':\ (z',z_n)\ \in S_\ell\}$ is contained in 
$$
E_{z_n}':=(\mathbb S_{\ell\de}^{n-2}+O(\de))\cap\ (\mathbb S^{n-2}_{x_n+z_n}-x'+O(\de)).
$$
For $|z_n|\in I_k:=[(k-1)\de,k\de]$ with $1\le k\le \ell+5$,
we have
\begin{equation}\label{e007}
	|E_{z_n}'|\les \ell^{\tf{n-2}2}(\ell+10-k)^{\tf{n-4}2}\de^{n-1}
\end{equation}
Let 
$$
A_{z_n}:=(\mathbb S_{\ell\de}^{n-2}+O(\de))\cap\ (\mathbb S^{n-2}_{x_n+z_n}-x').
$$
As $x_n+z_n\sim 1$,
\begin{equation}\label{e006}
|E_{z_n}'|\les |A_{z_n}|\de.
\end{equation}
Let $\al_{z_n,\ell}$ be the smallest angle of the triangle with side lengths $x_n+z_n,\ x_n,$ and $\ell\de$, i.e. the opposite angle of the side with length $\ell\de$. 
Then
$$
|A_{z_n}|\les [ \al_{z_n,\ell+1}]^{n-2}-[ \al_{z_n,\ell}]^{n-2}
\les \ell^{\tf{n-2}2}(\ell+10-k)^{\tf{n-4}2}\de^{n-2},
$$
which implies \eqref{e007} recalling \eqref{e006}.

It follows from \eqref{e007} that
\begin{align*}
&\int_{S_\ell}|z'|^{\al-(n-1)}|z_n|^{\be-1}dz\\
\les&\sum_{1\le k\le \ell+5}\int_{|z_n|\in I_k}\int_{\{z':\ (z',z_n)\ \in S_\ell\}}|z'|^{\al-(n-1)}dz'|z_n|^{\be-1}dz_n\\
\les &\sum_{1\le k\le \ell+1}\ell^{\tf{n-2}2}(\ell+10-k)^{\tf{n-4}2}\de^{n-1}(\ell\de)^{\al-(n-1)}(k\de)^{\be-1}\de\\
\sim &\de^{\al+\be}\ell^{\al-(n-1)}\ell^{\tf{n-2}2}\sum_{1\le k\le \ell+1}(\ell+10-k)^{\tf{n-4}2}k^{\be-1}\\
\les&\de^{\al+\be}\ell^{\al+\be-2}.
\end{align*}
Summing over $\ell$, we obtain the control
\begin{align*}
	\sum_\ell \int_{S_\ell}|z'|^{\al-(n-1)}|z_n|^{\be-1}dz
	\les &\de^{\al+\be}+\sum_{1\le \ell\les\de^{-1}}\de^{\al+\be}\ell^{\al+\be-2}+\de\\
	\les &\left\{\begin{array}{ll} 
		\de^{\al+\be},&\q \al+\be<1\\
		\de\log\tf1\de,&\q \al+\be=1\\
		\de,& \q \al+\be>1.
	\end{array}\right.
\end{align*}
This estimate 
finishes the proof.
\epf

\section{Weighted estimates}\label{t050}
The main goal of this section is to prove the following result via trace inequalities obtained in Section~\ref{040}.

\begin{thm}\label{113}
For $\al\in[0,n-1)$, $\be\in[0,1)$, and $0<\de\ll1$, let 
$$
A_\de(\al, \be)=\de^{\tf{1-(\al+\be)}2}C_\de(\al,\be)^{1/2}\sim \left\{\begin{array}{ll} 
	\de^{1/2}, &\q \al+\be<1\\
	\de^{1/2} (\log\tf1\de)^{1/2}, &\q \al+\be=1\\
	\de^{(2-\al-\be)/2}, & \q\al+\be>1.
\end{array}\right.$$
Then
\begin{equation}\label{e012}
\|G_0(f)\|_{L^2(w_{\al,\be})}\le C A_\de(\al,\be)\|f\|_{L^2(w_{\al,\be})}.
\end{equation}

\end{thm}

We recall that 
$w_{\al,\be}(x)=|x'|^{-\al}|x_n|^{-\be}$,
$$
m_\de(\xi)=\mu_\de(\xi'/\xi_n)\psi(2^{-1}\xi_n)
$$
with $\mu_\de$ satisfying \eqref{e008} and \eqref{e046},
and 
$$
G_0(f)(x)=\Big(\int_{1}^{2}\Big|\big[m_\de(t^{-1}\xi',\xi_n)\wh f(\xi)\big]^\vee(x)\Big|^2\f{dt}t\Big)^{1/2}.
$$

To study $G_0$, we decompose $K_\de:=(m_\de)^\vee$ dyadically. More precisely, we take a smooth function $\Psi$ supported in the unit annulus such that
$$
1=\sum_{j\ge j_0}\Psi_j,
$$
where 
$\Psi_j=\Psi(2^{-j}\cdot)$ is supported in $B(0,2^{j+1})\setminus B(0,2^{j-1})$ for $j\ge j_0+1$, and $\Psi_{j_0}=1-\sum_{j\ge j_0+1}\Psi_j$ is supported in $B(0,2^{j_0})$ with $2^{j_0}\sim\de^{-1}$. Correspondingly we decompose
$K_\de=\sum_{j\ge j_0}K_j$ with $K_j:= K_\de \Psi_j$.

As $m_\de$ is supported in $\Ga_\de$ and $\wh\Psi_j$ is essentially supported in a ball of radius $2^{-j}\les \de$, we see that
$\wh K_j=m_\de*\wh\Psi_j$ is essentially supported in ${\Ga_\de}$, which we describe precisely below.

\begin{prop}\label{06012'}
Let $j\ge j_0$. For any $M\in \mathbb N$, there exists a constant $C_M>0$ such that the following statements hold.

(i) For any $\xi\in\rn$, we have
\begin{equation}\label{e014}
|\wh{K_j}(\xi)|\le C_M(2^j\de)^{-M}.
\end{equation}

(ii) For dist $(\xi,\Ga)\sim 2^\ell\de$ with $\ell\ge 1$, we have
\begin{equation}\label{e010}
|\wh{K_j}(\xi)|\le C_M 2^{-\ell M}(2^j\de)^{-M}.
\end{equation}

(iii) For dist $(\xi,\Ga)\ge \tf12$, 
\begin{equation}\label{e009}
|\wh{K_j}(\xi)|\le C_M 2^{-j M}(1+|\xi|)^{-M}.
\end{equation}
\end{prop}

\bpf
(i) For $j\ge j_0+1$, we may rewrite $\wh K_j$ as 
$$
\int m_\de(\xi-2^{-j}\eta)\wh\Psi(\eta)d\eta.
$$
As one can easily show $\p^\al m_\de(\xi)\le\de^{-|\al|}$ from \eqref{e008} and \eqref{e046},
applying Taylor's expansion to $m_\de(\xi-2^{-j}\cdot)$, and using    the vanishing moments of $\Psi$, we obtain the claimed estimate.
When $j=j_0$, we have
$$
|\wh K_{j_0}(\xi)|\le C\|m_\de\|_{L^\nf}\|\wh \Psi_{j_0}\|_{L^1}\le C.
$$

(ii) When $\xi$ is $2^\ell\de$ away from $\Ga$, we obtain from the condition  $\xi-2^{-j}\eta\in \Ga_\de$ that $2^{-j}|\eta|\ges 2^\ell\de$. This implies  
\begin{align*}
\left|\int m_\de(\xi-2^{-j}\eta)\wh\Psi(\eta)d\eta\right|\le &C_{M,\Psi}\int_{|\eta|\ges 2^{j+\ell}\de}(1+|\eta|)^{-2M}d\eta\\
\le & C_{M,\Psi} (1+2^{j+\ell}\de)^{-M}\\
\le & C_{M,\Psi} 2^{-\ell M}(2^{j}\de)^{-M}.
\end{align*}

(iii) \eqref{e009} could be proved similarly via the observation that, for $\xi-2^{-j}\eta\in\Ga_\de$, we have
$$
2^{-j}|\eta|\ge \tf1{10}(1+|\xi|)
$$ 
due to the assumption dist $(\xi,\Ga)\ge \tf12$.
\epf

We define $\wh{K_{\de,t}}(\xi)=m_\de(t^{-1}\xi',\xi_n)$, which implies that $K_{\de,t}=t^{n-1}K_\de(tx',x_n)$. Let $K_{j,t}(x)=t^{n-1}K_j(tx',x_n)$, 
 then
\begin{equation}\label{e05211}
K_{\de,t}(x)=\sum_{j\ge j_0}K_{j,t}(x).
\end{equation}
Obviously 
\begin{equation}\label{e01291'}
G_0(f)(x)\le\sum_{j\ge j_0}\Big(\int_{1}^{2}|(K_{j,t}*f)(x)|^2\f{dt}t\Big)^{1/2}.
\end{equation}

We remark that $K_{j,t}$ is supported in a ball of radius $c_1 2^j$ with $c_1\sim 1$ as $t\in[1,2]$. Therefore we can localize $f$ by decomposing it into pieces supported in cubes of length $c_1 2^j$. 
For $i\in\bbz^n$, we define 
$  Q_i=Q(c_12^ji, c_12^j)$,  then $\{  Q_i\}_{i\in\bbz^n}$ is a partition of $\rn$.
Define $f_i=f\chi_{Q_i}$, then \begin{equation}\label{e013}
f=\sum_{i\in\bbz^n}f_i.
\end{equation}

To prove \eqref{e012}, it follows from \eqref{e01291'} and \eqref{e013} that we should estimate
$$
\int_{1}^{2}\int_{\rn}|K_{j,t}*f_i(x)|^2w_{\al,\be}(x)dx\f{dt}t,
$$
which relies heavily on the position of $i$. Therefore we decompose $\bbz^n$ into four subsets:
\begin{eqnarray*}
E_1:=\{i\in\bbz^n:\ |i'|\ge 10n, |i_n|\ge 10n\}\\
E_2:=\{i\in\bbz^n:\ |i'|\ge 10n, |i_n|\le 10n\}\\
E_3:=\{i\in\bbz^n:\ |i'|\le 10n, |i_n|\le 10n\}\\
E_4:=\{i\in\bbz^n:\ |i'|\le 10n, |i_n|\ge 10n\},
\end{eqnarray*}
and discuss them separately.

\medskip

\noindent{\bf Case 1: $i\in E_1$}

\medskip

To handle this case, we need the following lemma obtained from the Plancherel identity.

\begin{lm}\label{107}

Suppose 
$0< A\le w(x)\le B$ for 
all $x\in 3Q_i$, where $  3Q_i=Q(c_12^ji, 3c_12^j)$. 
Then
\begin{equation}\label{e109}
\Big\|\Big(\int_{1}^{2}|K_{j,t}*f_i(x)|^2\f{dt}t\Big)^{1/2}\Big\|_{L^2(w)}
\les
(\tf BA)^{1/2}\de^{1/2}(2^j\de)^{-M}\|f_i\|_{L^2(w)}.
\end{equation}
\end{lm}

\bpf
Noticing that $K_{j,t}*f_i$ is supported in $3Q_i$, the left hand side of \eqref{e109}  is bounded by 
\begin{equation}\label{e01274'}
B^{1/2}\Big(\int_{1}^{2}\int_{\rn}|K_{j,t}*f_i(x)|^2dx\f{dt}t\Big)^{1/2}.
\end{equation}

As $\wh {K_{j,t}}$ is essentially supported{\footnote{We present here an intuitive argument, while a strict one is provided in the appendix.}} in 
$$\Ga_{\de,t}:=\{\xi\in\rn:\ (t^{-1}\xi',\xi_n)\in\Ga_\de\},
$$
 we may decompose the range of $t$ as 
 $$[1,2]=\cup_{\be=\de^{-1}}^{2\de^{-1}}I_\be
 $$ 
 with $I_\be=[\be  \de,(\be+1) \de]$. Moreover, we define $\wh {P_\be g}$ as the restriction of $\wh g$ to the set
$$
\{\xi\in \rn: \tf{|\xi'|}{\xi_n}\in[\be  \de,(\be+1) \de], \xi_n\in[\tf12,4]\}.
$$
It follows from Proposition~\ref{06012'} that, for $t\in I_\be$, the essential support of $\wh{K_{j,t}}$ and the support of $\wh {P_{\be'}g}$ are  disjoint when $|\be-\be'|\ge 10$.
This implies that $K_{j,t}*P_{\be'}f_i$ is essentially $0$ when $t\in I_\be$ with $|\be-\be'|\ge 10$.
Hence 
\begin{align*}
&\int_{1}^{2}\int_{\rn}|K_{j,t}*f_i(x)|^2dx\f{dt}t\\
=&\sum_\be\int_{I_\be}\int_{\rn}|K_{j,t}*f_i(x)|^2dx\f{dt}t\\
\sim &\sum_\be\int_{I_\be}\int \Big|\sum_{\be'=\be-10}^{\be+10}(K_{j,t}*(P_{\be'}f_i))(x)\Big|^2dx\f{dt}t\\
\les &\sum_{\be}\sum_{\be'=\be-10}^{\be+10}\int_{I_\be}\int |(K_{j,t}*(P_{\be'}f_i))(x)|^2dx\f{dt}t\\
\sim &\sum_{\be}\sum_{\be'=\be-10}^{\be+10}\int_{I_\be}\int |(\wh{K_{j,t}}\wh{P_{\be'}f_i})(\xi)|^2d\xi \f{dt}t\\
\les &(2^j\de)^{-M}\sum_{\be}\sum_{\be'=\be-10}^{\be+10}\int_{I_\be}\int |(\wh{P_{\be'}f_i})(\xi)|^2d\xi \f{dt}t\\
\les &\de (2^j\de)^{-M}\sum_{\be}\sum_{\be'=\be-10}^{\be+10} \int |(\wh{P_{\be'}f_i})(\xi)|^2d\xi \\
\les &\de (2^j\de)^{-M}\|f_i\|_{L^2}^2\\
\les &A^{-1}\de (2^j\de)^{-M}\|f_i\|_{L^2(w)}^2,
\end{align*}
where we use \eqref{e014}.
This combined with \eqref{e01274'} yields \eqref{e109}.
\epf

For $i\in E_1$ and $w(x)=w_{\al,\be}(x)$,  we have $A\sim B$. Applying Lemma~\ref{107}, we obtain 
\begin{equation}\label{e110}
\Big(\int_{1}^{2}\int_{\rn}|K_{j,t}*f_i(x)|^2w_{\al,\be}(x)dx\f{dt}t\Big)^{1/2}\le C\de^{1/2}(2^j\de)^{-M}\|f_i\|_{L^2(w_{\al,\be})}.
\end{equation}

\medskip

\noindent{\bf Case 2: $i\in E_2$}

\medskip

In the remaining three cases, as $w_{\al,\be}$ may not be bounded, Lemma~\ref{107} is not applicable. We need the following estimate to make use of
 trace inequalities obtained in Section~\ref{040}.

\begin{lm}\label{01273'}

Suppose that $w_1(x)\ge 0$ satisfies
\begin{equation}
 w_1(x)\ge A>0  \qq \forall x\in Q_i, 
\end{equation}
\begin{equation}\label{e203}
w_1(x)\sim w_1(sx',x_n)\quad \text{for } s\sim 1, 
\end{equation}
and 
\begin{equation}\label{e01275'}
\int_{\Ga_\de}|\wh g(\xi)|^2 d\xi\le C_\de(w_1)\int_{\rn}|g(x)|^2 w_1^{-1}(x)dx .
\end{equation}

Then 
$$
\Big(\int_{1}^{2}\int_{\rn}|K_{j,t}*f_i(x)|^2w_1(x)dx\f{dt}t\Big)^{1/2}
\le 
CA^{-1/2}\de^{1/2}(2^j\de)^{-M}C_\de(w_1)^{1/2}\|f_i\|_{L^2(w_1)}.
$$

\end{lm}

\bpf
As in the proof of Lemma~\ref{107}, we write $[1,2]=\cup_{\be=\de^{-1}}^{2\de^{-1}}I_\be$ with $I_\be=[\be \de,(\be+1)\de]$. Then
\begin{align}
&\int_{1}^{2}\int_{\rn}|K_{j,t}*f_i(x)|^2w_1(x)dx\f{dt}t\notag\\
\le &C\sum_\be\int_{I_\be}\int_{\rn} |\sum_{\be'=\be-10}^{\be+10}(K_{j,t}*(P_{\be'}f_i))(x)|^2w_1(x)dx\f{dt}t\notag\\
\le&C \sum_\be\sum_{\be'=\be-10}^{\be+10}\int_{I_\be}\int_{\rn} |(K_{j,t}*(P_{\be'}f_i))(x)|^2w_1(x)dx\f{dt}t.\label{e01276'}
\end{align}

Applying \eqref{e01275'} and Proposition~\ref{06012'}, we obtain essentially
\begin{align*}
\int_{\rn}|K_{j,t}*g(x)|^2dx=&\int_{\rn} |\wh{K_{j,t}}(\xi)\wh g(\xi)|^2d\xi\\
\les&(2^j\de)^{-M}\int_{\Ga_{\de,t}} |\wh g(\xi)|^2 d\xi\\
\les& (2^j\de)^{-M} C_\de(w_1)\int_{\rn} |g(x)|^2w_1^{-1}(x)dx,
\end{align*}
where in the last step we use \eqref{e203}.
By duality
$$
\int_{\rn}|K_{j,t}*g(x)|^2w_1(x)dx\les (2^j\de)^{-M} C_\de(w_1)\int_{\rn} |g(x)|^2dx.
$$
Hence \eqref{e01276'} is bounded by 
\begin{align*}
&(2^j\de)^{-M}C_\de(w_1)\sum_\be\sum_{\be'=\be-10}^{\be+10}\int_{I_\be}\int_{\rn} |P_{\be'}f_i(x)|^2dx\f{dt}t\\
\les&(2^j\de)^{-M}\de C_\de(w_1)\|f_i\|_{L^2}^2\\
\les&(2^j\de)^{-M}\de C_\de(w_1)A^{-1}\|f_i\|_{L^2(w_1)}^2.
\end{align*}
This completes the proof.

\epf

We observe that when $x\in Q_i$ with $i\in E_2$, we have  $|x_n|\les 2^j$ and $|x'|^{-\al}\sim |x_i'|^{-\al}$, where $x_i=c_12^ji$ is the center of $Q_i$.
Take $w_1(x)=|x_n|^{-\be}$, and accordingly $A=2^{-j\be}$, then $C_\de(w_1)=\de^\be$ by Proposition~\ref{109'} with $\al=0$. Applying Lemma~\ref{01273'} with $M_1\ge M+\be$, we obtain
\begin{align*}
&\int_{1}^{2}\int_{\rn}|K_{j,t}*f_i(x)|^2w_{\al,\be}(x)dx\f{dt}t\\
\les&|x_i'|^{-\al}\int_{1}^{2}\int_\rn|K_{j,t}*f_i(x)|^2w_{1}(x)dx\f{dt}t\\
\les&|x_i'|^{-\al}\de 2^{j\be}(2^j\de)^{-M_1}\de^{\be}\int_\rn|f_i(x)|^2w_1(x)dx\\
\les&\de(2^j\de)^{-M}\int_\rn|f_i(x)|^2w_{\al,\be}(x)dx,
\end{align*}
as $w_{\al,\be}(x)\sim |x_i'|^{-\al}w_1(x)$ in $Q_i$.
In summary, we have
\begin{equation}\label{e111}
\Big(\int_{1}^{2}\int_{\rn}|K_{j,t}*f_i(x)|^2w_{\al,\be}(x)dx\f{dt}t\Big)^{1/2}
\le 
C\de^{1/2}(2^j\de)^{-M}\|f_i\|_{L^2(w_{\al,\be})}.
\end{equation}

\medskip

\noindent{\bf Case 3: $i\in E_3$}

\medskip

In this case, to apply Lemma~\ref{01273'}, we take $w_1(x)=w_{\al,\be}(x)$ and $A=2^{-j(\al+\be)}$ such that $A\le w_1(x)$ in $Q_i$. We recall that by Proposition~\ref{109'}
$$
C_\de(w_{\al,\be})\les \left\{\begin{array}{ll} 
\de^{\al+\be},&\q \al+\be<1\\
\de\log\tf1\de,&\q \al+\be=1\\
\de,& \q \al+\be>1,
\end{array}\right.$$
therefore it follows from Lemma~\ref{01273'} that 
\begin{align}\label{e112}
&\Big(\int_{1}^{2}\int_{\rn}|K_{j,t}*f_i(x)|^2w_{\al,\be}(x)dx\f{dt}t\Big)^{1/2}\notag\\
\le &
C2^{j(\al+\be)/2}\de^{1/2}(2^j\de)^{-M}C_\de(w_{\al,\be})^{1/2}\|f_i\|_{L^2(w_{\al,\be})}\notag\\
\les&\|f_i\|_{L^2(w_{\al,\be})} \left\{\begin{array}{ll} 
\de^{1/2}(2^j\de)^{-M},&\q \al+\be<1\\
\de^{1/2}(\log\tf1\de)^{1/2}(2^j\de)^{-M},&\q \al+\be=1\\
\de^{1-\tf{\al+\be}2}(2^j\de)^{-M},& \q \al+\be>1.
\end{array}\right.
\end{align}

\medskip

\noindent{\bf Case 4: $i\in E_4$}

\medskip

We recall that $E_4:=\{i\in\bbz^n:\ |i'|\le 10n, |i_n|\ge 10n\}$, therefore  in $Q_i$ with $i\in E_4$ we have $|x'|\les 2^j$ and $|x_n|^{-\be}\sim (2^j|i_n|)^{-\be}$. Taking $w_1(x)=|x'|^{-\al}=w_{\al,0}$ and $A=2^{-j\al}$, we can apply
 Lemma~\ref{01273'} to obtain
\begin{align*}
&\int_{1}^{2}\int_{\rn}|K_{j,t}*f_i(x)|^2w_{\al,\be}(x)dx\f{dt}t\\
\les&(2^j|i_n|)^{-\be}\int_{1}^{2}\int_\rn|K_{j,t}*f_i(x)|^2w_{1}(x)dx\f{dt}t\\
\les&(2^j|i_n|)^{-\be}\de 2^{j\al}(2^{j}\de)^{-M}C_\de(w_1)\int_\rn|f_i(x)|^2w_1(x)dx\\
\les&\de 2^{j\al}(2^{j}\de)^{-M}C_\de(\al,0)\int_\rn|f_i(x)|^2w_{\al,\be}(x)dx.
\end{align*}
By Proposition~\ref{109'} with $\be=0$ we obtain
\begin{align}\label{e113}
&\Big(\int_{1}^{2}\int_{\rn}|K_{j,t}*f_i(x)|^2w_{\al,\be}(x)dx\f{dt}t\Big)^{1/2}\notag\\
\les &\|f_i\|_{L^2(w_{\al,\be})} \left\{\begin{array}{ll} 
\de^{1/2}(2^j\de)^{-M},&\q0<\al<1\\
\de^{1/2}(\log\tf1\de)^{1/2}(2^j\de)^{-M},&\q \al=1\\
\de^{1-\tf{\al}2}(2^j\de)^{-M},& \q 1<\al<n-1.
\end{array}\right.
\end{align}

We are now ready to prove Theorem~\ref{113}.

\bpf[Proof of Theorem~\ref{113}]

When $\al=\be=0$, \eqref{e012} is verified by Proposition~\ref{117}.
We discuss below the case $(\al,\be)\neq (0,0)$.

Recalling that $f_i$ is supported in $Q_i$ and $K_{j,t}*f_i$ is supported in $3Q_i$, we have
$K_{j,t}*f=\sum_{i}(K_{j,t}*f_i)\chi_{3Q_i}\les (\sum_i |K_{j,t}*f_i|^2)^{1/2}$ as $\{3Q_i\}_{i\in\bbz^n}$ are finitely overlapping. 
It follows from the definition of $f_i$ that $|\sum_i f_i|^2=\sum_i |f_i|^2$.
Hence 
\begin{align*}
&\int_{1}^{2}\int_\rn|K_{j,t}*f|^2 |x'|^{-\al}|x_n|^{-\be}dx\f{dt}t\\
\les&\sum_i \int_{1}^{2}\int_\rn|K_{j,t}*f_i|^2  |x'|^{-\al}|x_n|^{-\be}dx\f{dt}t\\
\sim&\Big(\sum_{i\in E_{1}}+\sum_{i\in E_{2}}+\sum_{i\in E_{3}}+\sum_{i\in E_{4}}\Big)
\int_{1}^{2}\int_\rn|K_{j,t}*f_i|^2 |x'|^{-\al}|x_n|^{-\be}dx\f{dt}t.
\end{align*}

When $\al+\be>1$, by \eqref{e112}, we have, for $i\in E_3$,
$$
\int_{1}^{2}\int_\rn|K_{j,t}*f_i|^2 |x'|^{-\al}|x_n|^{-\be}dx\f{dt}t 
\le \de^{2-\al-\be}(2^j\de)^{-M}\|f_i\|^2_{L^2(w_{\al,\be})}.
$$
Moreover, noticing that $1\le C\de^{-\be}$, the same bound holds 
for all $\al\in(0,n-1)$ when $i\in E_4$ by \eqref{e113}.
These estimates and \eqref{e110}, \eqref{e111} yield
\begin{align*}
&\int_{1}^{2}\int_\rn|K_{j,t}*f|^2 |x'|^{-\al}|x_n|^{-\be}dx\f{dt}t\\
\les&\big[\de(2^j\de)^{-M}+\de^{2-\al-\be}(2^j\de)^{-M}\big]
(\sum_{i\in E_{1}}+\sum_{i\in E_{2}}+\sum_{i\in E_{3}}+\sum_{i\in E_{4}}) \int_\rn|f_i|^2 |x'|^{-\al}|x_n|^{-\be}dx \\
\les &\big[\de^{2-\al-\be}(2^j\de)^{-M}\big]  
 \int_\rn|f|^2  |x'|^{-\al}|x_n|^{-\be}dx .
\end{align*}
This combined with \eqref{e01291'}
implies that 
$$
\|G_0(f)\|_{L^2(w_{\al,\be})}\les\de^{(2-\al-\be)/2}\|f\|_{L^2(w_{\al,\be})}
$$
since $j\ge j_0$ with $2^{j_0}\sim\de^{-1}$.

The case when $\al+\be\le1$ can be proved similarly.

\epf

\section{Proof of main results}\label{t083}

It is necessary to show that $L$ is well defined in  $L^2(w_{\al,\be})$ and $L^p(\rn)$.

\begin{lm}\label{t012}
Let $\al\in[0,n-1)$, $\be\in [0,1)$, and $p\in (1,\nf)$.
The linear operator $L$ initially defined on Schwartz functions is bounded on $L^2(w_{\al,\be})$ and bounded on $L^p(\rn)$.

\end{lm}

\bpf
We recall that 
\begin{align*}
L(f)(x)=&\int_\rn \psi(2^{-1}\xi_n)\wh f(\xi)e^{2\pi ix\cdot\xi}d\xi\\
=&\int_{\bbr^{n-1}}\int_\bbr f(x',y_n)2\psi^\vee(2(x_n-y_n))dy_ndx'.
\end{align*}
The $L^p$ boundedness of $L$ follows from Minkowski's inequality directly as $\psi^\vee\in L^1(\bbr)$.

Observing that $|g*\psi^\vee(x_n)|\le CM(g)(x_n)$, which is bounded on $L^2(|x_n|^{-\be})$for $\be\in(-1,1)$, the boundedness of $L$ on $L^2(w_{\al,\be})$ follows.

\epf

\bpf[Proof of Theorem~\ref{06011}]

When $\al+\be>1$, we obtain from Theorem~\ref{113} and Proposition~\ref{101} that
$$
\|{\mathbf G}_\de(f)\|_{L^2(w_{\al,\be})}\les\de^{(2-\al-\be)/2}\|f\|_{L^2(w_{\al,\be})}.
$$
It follows from this, Proposition~\ref{111}, and \eqref{e01272'} that 
$$
\|T_*^\la(f)\|_{L^2(w_{\al,\be})}
\les\sum_{\ga\ge0}2^{-\ga\la}2^{\ga/2}2^{-\ga(2-\al-\be)/2}\|f\|_{L^2(w_{\al,\be})}
\les C\|f\|_{L^2(w_{\al,\be})}
$$
as $\la>\tf{\al+\be-1}2$.

The case when $\al+\be\le 1$ can be proved similarly with the help of Theorem~\ref{113}.

\medskip

We have finished the proof of \eqref{e033}.
Now we turn to the proof of the second statement, which follows from
the boundedness of $T^\la_*$ on $L^2(\rn,w_{\al,\be})$ and the following two lemmas.

\begin{lm}\label{06021}
For any measurable set $E$, $|E|=0$ is equivalent to $w_{\al,\be}(E)=0$.
\end{lm}

This lemma is trivial.

\begin{lm}[{\cite[Exercise 7.4.1]{Grafakos2014b}}]
Let $\al\in[0,n-1)$ and $\be\in [0,1)$.
The Schwartz functions are dense in $L^2(w_{\al,\be})$.

\end{lm}

By these lemmas, Lemma~\ref{t012} and \eqref{e033}, we  apply Stein's maximal principle (\cite{Stein1961a}, see also \cite[Theorem 2.1.14.]{Grafakos2014b}) to obtain that
$$
\lim_{t\to\nf}T^\la_t(f)(x)=L(f)(x) \qq a.e.
$$
whenever $f\in L^2(w_{\al,\be})$.

\epf

Next we prove Theorem~\ref{01276}.

\bpf[Proof of Theorem~\ref{01276}]

For a given $\la>0$, we fix $p\in[2,\f{2n}{n-2\la-1})$.
For any $f\in L^p(\rn)$, by Proposition~\ref{116} we can write $f=\sum_{i=1}^4f_i$ with $f\in L^2(w_i)$. 
In the decomposition, for any $\ep>0$, we can select $\al_i,\be_i>0$ appropriately such that $\al_i+\be_i<n(1-\tf2p)+\ep$, $i=1,2,3,4$. In particular, we have
$\la>(\al_i+\be_i-1)/2$ for $i=1,2,3,4$.

By Theorem~\ref{06011} we obtain 
$$
\lim_{t\to\nf}T^\la_t(f_i)(x)=L(f_i)(x) \qq a.e.
$$
for $i=1,2,3,4$. This complete the proof of Theorem~\ref{01276} by summing over $i$.
\epf

Finally we prove Theorem~\ref{t013}.

\bpf[Proof of Theorem~\ref{t013}]
Fix $k_1$ and $k_2$ in $\bbz$ such that $2^{k_1+10}<a<b<2^{k_2-10}$.

We define 
$$
L_k(f)(x)=\int_{\rn} \psi(2^{-k-1}\xi_n)\wh f(\xi)e^{2\pi ix\cdot\xi} d\xi.
$$
Then for $f\in L^p(\rn)$ satisfying  \eqref{e036}, $L_k$ is well defined, and
\begin{equation}\label{e037}
f=\sum_{k_1\le k\le k_2} L_k(f).
\end{equation}
Let 
$$
 T_{k,t}^\la(f)(x):=\int_{\rn} (1-\tf{|\xi'|^2}{\xi_n^2})^\la_+\psi(2^{-k-1}\xi_n)\wh f(\xi) e^{2\pi ix\cdot\xi}d\xi,
$$
then a  dilation argument shows that 
\begin{equation}\label{e038}
\lim_{t\to\nf} T_{k,t}^\la(f)= L_k(f)\qq a.e.
\end{equation}
Noticing that, for $f\in L^p(\rn)$ satisfying  \eqref{e036},
$$
\tilde T^\la_t(f)=\sum_{k_1\le k\le k_2}T_{k,t}^\la(f),
$$
the conclusion follows by \eqref{e037} and \eqref{e038}.

\epf

\section{Remarks}\label{t084}

\subsection{}
By a standard duality argument, we  obtain from Lemma~\ref{103} the following orthogonality result.

\begin{lm}\label{104}

Let $k,L\in\bbz$, $\al\in(-(n-1),n-1)$ and $\be\in(-1,1)$. There exists a constant $C$ such that, for any sequence of functions $\{H_{k,L}\}$ satisfying 
$$
\text{supp }\wh{H_{k,L}}\sset \{\xi\in\rn:\ |\xi'|\sim 2^{L+k},\xi_n\sim 2^L\},
$$
we have
$$
\int_\rn|\sum_{L\in\bbz}H_{k,L}(x)|^2w_{\al,\be}(x)dx
\le C \sum_{L\in\bbz}\int_\rn |H_{k,L}(x)|^2w_{\al,\be}(x)dx.
$$

\end{lm}

Let 
$$
\widetilde
{\mathbf G}_\de(f)(x):=\left(\int_0^\nf\Big|\int_{\rn}\mu_\de(\xi'/t\xi_n)\wh f(\xi)e^{2\pi ix\cdot\xi}d\xi\Big|^2\f{dt}t\right)^{1/2}.
$$
With the help of Lemma~\ref{104}, we can actually obtain from \eqref{e048} the estimate
\begin{equation}\label{e060}
\|\widetilde {\mathbf G}_\de(f)\|_{L^2(w_{\al,\be})}\les A\|f\|_{L^2(w_{\al,\be})}
\end{equation}
by imitating the proof of Proposition~\ref{101}.

By \eqref{e060}, to obtain the (weighted) boundedness of 
$$
\sup_{t>0}\Big|\int_{\rn} (1-\tf{|\xi'|^2}{t^2\xi_n^2})^\la_+\wh f(\xi) e^{2\pi ix\cdot\xi}d\xi\Big|,
$$
it remains to study the maximal function defined by 
\begin{equation}\label{e086}
\widetilde M_\vp(f)(x):=\sup_{t>0}\Big|[\vp(\xi'/t\xi_n)\wh f(\xi)]^\vee\Big|,
\end{equation}
where $\vp\in C_0^\nf(\bbr^{n-1})$ is supported in $B(0,1)$ with $\vp(0)=1$. We conjecture that 
$\widetilde M_\vp$ is bounded on $L^2(\rn)$. But we cannot prove it. By a standard argument using square functions, we can  assume further that $\vp(\xi')=1$ for $|\xi'|\le 1/2$.

\subsection{}
Let 
$$
K_\la(x)=\int_\rn (1-\tf{|\xi'|^2}{\xi_n^2})^\la_+\psi(\xi_n) e^{2\pi ix\cdot\xi }d\xi.
$$
It has been proved in \cite{Heo2010} that 
$$
K^\la\in L^1(\rn)
$$
for $\la>\tf {n-2}2$, which implies that 
\begin{equation}
\|T^\la(f)\|_{L^p(\rn)}\le C\|f\|_{L^p(\rn)}
\end{equation}
for $1<p<\nf$.
It is natural to conjecture that, when $\la>\tf {n-2}2$,
\begin{equation}\label{e088}
\|T^\la_*(f)\|_{L^p(\rn)}\le C\|f\|_{L^p(\rn)},\q 1<p<\nf.
\end{equation}
This can be verified for $\la>\tf n2$

\begin{prop}\label{t085}
Let $\la>\tf n2$. Then
\begin{equation*}
\|T^\la_*(f)\|_{L^p(\rn)}\le C\|f\|_{L^p(\rn)},\q 1<p<\nf.
\end{equation*}
\end{prop}

This result follows from the pointwise estimate
\begin{equation}\label{e089}
|K_\la(x)|\le C(1+|x|)^{-\tf n2-\la},
\end{equation}
obtained by \cite{Heo2010}.
We remark that this bound is worse than the corresponding pointwise estimate of $\big[(1-|\xi|^2)_+\big]^\vee$.

We are curious if we can fill the gap between Proposition~\ref{t085} and the conjectured range $\la>\tf {n-2}2$. This kind of gap  does appear in the classical Bochner-Riesz problem.

\section{Appendix}

We present a strict proof of Lemma~\ref{107} in this section, while a strict proof of Lemma~\ref{01273'} could be obtained similarly.

Let us recall Lemma~\ref{107} first.

\begin{lm}

Suppose 
$ 0<A\le w(x)\le B$ for 
all $x\in 3Q_i$, where $  3Q_i=Q(c_12^ji, 3c_12^j)$. 
Then
$$
\Big\|\Big(\int_{1}^{2}|K_{j,t}*f_i(x)|^2\f{dt}t\Big)^{1/2}\Big\|_{L^2(w)}
\les
(\tf BA)^{1/2}\de^{1/2}(2^j\de)^{-M}\|f_i\|_{L^2(w)}.
$$

\end{lm}

\bpf

Let $\ell_0$ be an integer such that $2^{\ell_0}\sim (10\de)^{-1}$.
We need to decompose $\wh K_j$ further such that
$$
\wh K_j(\xi)=\sum_{10\le\ell\le\ell_0}\wh K_j(\xi)\wh\vp_\ell(\xi)+\wh K_j(\xi)\wh\vp_\nf(\xi)=:\sum_{\ell}\wh{K_{j,\ell}}(\xi),
$$
where $\wh\vp_{10}(\xi)$ is supported in $\{\xi:\ dist(\xi,\Ga)\le  2^{20} \de\}$,  $\wh{\vp_{\ell}}$ is supported in $\{\xi:\ dist (\xi,\Ga)\sim 2^\ell\de\}$ for $11\le\ell\le \ell_0$, and 
$$
\wh\vp_\nf(\xi)=1-\sum_{10\le\ell\le\ell_0}\wh\vp_\ell(\xi).
$$
Correspondingly we define $\wh {K_{j,\ell,t}}(\xi)=\wh {K_{j,\ell}}(t^{-1}\xi',\xi_n)$, and we have
$$
K_{j,t}=\sum_{\ell\ge 10}K_{j,\ell,t}.
$$

Using the above decomposition and that $ w(x)\le B$ we conrtol $(\int_1^{2}\int_{\rn}|K_{j,t}*f_i(x)|^2w(x)dx\f{dt}t)^{1/2}$ by 
\begin{equation}\label{e032}
B^{1/2}\sum_{\ell\ge 10}\Big(\int_{1}^{2}\int_{\rn}|K_{j,\ell,t}*f_i(x)|^2dx\f{dt}t\Big)^{1/2}.
\end{equation}
 We will consider two cases: $10\le \ell\le \ell_0$, and $\ell\ge \nf$.

\medskip

\noindent {\bf Case 1}: $10\le \ell\le\ell_0$.

\medskip

Since $\wh {K_{j,\ell,t}*f_i}(\xi)=\wh{K_{j,\ell}}(t^{-1}\xi',\xi_n)\wh f_i(\xi)$, in its support we have 
$$
|\xi'|=t\xi_n+[-2^{\ell+2}\de,2^{\ell+2}\de]
$$
and  $\xi_n\in[\tf12,4]$.
We decompose $[1,2]=\cup_{\be=\de^{-1}}^{2\de^{-1}} I_\be$ with $I_\be=[\be  \de,(\be+1) \de]$, and
 the support of $\wh {K_{j,\ell,t}}$ is 
$$
\{\xi\in\rn:\ |\xi'|\in\xi_n\be  \de+[-2^{\ell+2}\de,2^{\ell+2}\de],\ \xi_n\in[\tf12,4]\}.
$$
Define $\wh{P_\be g}$ as the restricition of $\wh g$ to the set 
$$
\{\xi\in\rn:\ |\xi'|\in[\xi_n\be  \de,\xi_n(\be+1)\de],\ \xi_n\in[\tf12,4]\}.
$$
Then for $t\in I_\be$, the support of $\wh {K_{j,\ell,t}}$ and the support of $\wh {P_{\be'}g}$ are disjoint if $|\be-\be'|\ge 10\cdot 2^\ell$. 
Hence
\begin{align}
&\int_1^2\int |K_{j,\ell,t}*f_i(x)|^2dx\f{dt}t\notag\\
=&\sum_{\be}\int_{I_\be}\int \big|\sum_{\be'=\be-10\cdot 2^\ell}^{\be+10\cdot 2^\ell}(K_{j,\ell,t}*P_{\be'}f_i)(x)\big|^2dx\f{dt}t\notag\\
\le& C2^\ell \sum_\be\ \sum_{\be'=\be-10\cdot 2^\ell}^{\be+10\cdot 2^\ell}\int_{I_\be}\int \big|(K_{j,\ell,t}*P_{\be'}f_i)(x)\big|^2dx\f{dt}t,\label{e031}
\end{align}
where in the last step we use H\"older's inequality.
By Proposition~\ref{06012'} (i) and (ii) we have
\begin{align*}
\int \big|(K_{j,\ell,t}*P_{\be'}f_i)(x)\big|^2dx=&\int \big|\wh{K_{j,\ell,t}}(\xi)\wh{P_{\be'}f_i}(\xi)\big|^2d\xi\\
\les& 2^{-\ell M}(2^j\de)^{-M}\int \big|\wh{P_{\be'}f_i}(\xi)\big|^2d\xi.
\end{align*}
Hence we control \eqref{e031} by 
\begin{align*}
&2^\ell 2^{-\ell M}(2^j\de)^{-M} \sum_\be\ \sum_{\be'=\be-10\cdot 2^\ell}^{\be+10\cdot 2^\ell}
\int_{I_\be}\int \big|\wh{P_{\be'}f_i}(\xi)\big|^2d\xi\f{dt}t\\
\les& 2^\ell 2^{-\ell M}(2^j\de)^{-M}2^\ell \de\sum_{\be'}\int \big|\wh{P_{\be'}f_i}(\xi)\big|^2d\xi\\
\les&2^{-\ell (M-2)}\de (2^j\de)^{-M}\|f_i\|_{L^2}^2.
\end{align*}

\noindent{\bf Case 2}: $\ell= \nf$

In this case we apply Proposition~\ref{06012'} (iii) directly to obtain the control
\begin{align*}
&\int_1^2\int |K_{j,\ell,t}*f_i(x)|^2dx\f{dt}t\\
\les&2^{-j M}\|f_i\|_{L^2}^2\\
\les& \de (2^j\de)^{-M}\|f_i\|_{L^2}^2,
\end{align*}
where in the last step we use that $\de\le 1$.

\medskip

By the estimates obtained in the above two cases, we can control \eqref{e032} by 
$$
B^{1/2}(\de(2^j\de)^{-M})^{1/2}\sum_{\ell\ge 10}2^{-\ell M}\|f_i\|_2
\les (B/A)^{1/2}\de^{1/2}(2^j\de)^{-M/2}\|f_i\|_{L^2(w)}
$$
since $w(x)\ge A$ in the support of $f_i$.
\epf

\end{document}